%% file: entropies-of-SFTs.arxiv.tex
\documentclass[oneside,english]{amsart}
\usepackage[latin1]{inputenc}
\usepackage{setspace}
\usepackage{amssymb}

\makeatletter

\providecommand{\LyX}{L\kern-.1667em\lower.25em\hbox{Y}\kern-.125emX\@}

 \theoremstyle{plain}    
 \newtheorem{thm}{Theorem}[section]
 \numberwithin{equation}{section} 
 \numberwithin{figure}{section} 
 \theoremstyle{remark}    
 \newtheorem*{acknowledgement*}{Acknowledgement} 
 \theoremstyle{plain}    
 \newtheorem{algorithm}[thm]{Algorithm} 
 \theoremstyle{plain}    
 \newtheorem{cor}[thm]{Corollary} 
 \theoremstyle{plain}    
 \newtheorem{lem}[thm]{Lemma} 
 \theoremstyle{definition}
  \newtheorem{problem}[thm]{Problem}
 \theoremstyle{plain}    
 \newtheorem{prop}[thm]{Proposition} 
 \theoremstyle{plain}    
 \newtheorem*{thm*}{Theorem} 

\usepackage{color}
\usepackage{epsfig}
\usepackage{extpic}

\newcommand{\hide}[1]{}
\DeclareMathOperator{\diam}{diam}

\usepackage{babel}
\makeatother
\begin{document}

\title{A characterization of the entropies of multidimensional shifts of
finite type}

\author{Michael Hochman and Tom Meyerovitch}

\begin{abstract}
We show that the values of entropies of multidimensional shifts of
finite type (SFTs) are characterized by a certain computation-theoretic
property: a real number $h\geq 0$ is the entropy of such an SFT if
and only if it is right recursively enumerable, i.e. there is a computable
sequence of rational numbers converging to $h$ from above. The same
characterization holds for the entropies of sofic shifts. On the other
hand, the entropy of an irreducible SFT is computable.
\end{abstract}

\email{mhochman@math.huji.ac.il, tomm@post.tau.ac.il}


\maketitle

\section{\label{sec:Introduction}Introduction}

A shift of finite type (SFT) is an ensemble of colorings of $\mathbb{Z}$
(a one-dimensional SFT) or $\mathbb{Z}^{d}$ for $d>1$ (a multidimensional
SFT) defined by local rules. SFTs are one of the fundamental objects
of study in symbolic dynamics, and their most significant invariant
is their (topological) entropy, which measures the asymptotic growth
of the number of legal colorings of finite regions (see section \ref{sec:Preliminaries}
for definitions). Besides having been been studied extensively from
a dynamical perspective as topological analogs of Markov chains \cite{P64,W73,W77},
SFTs appear naturally in a wide range of other disciplines. In information
theory, SFTs were used by Shannon as models for discrete communication
channels \cite{S48}, for which entropy describes the capacity; similarly,
SFTs model {}``two-dimensional'' channels \cite{FJ99}. SFTs have
been used to study the dynamics of geodesic flows and have played
an important role in the classification of the dynamics of Anosov
and Axiom A diffeomorphisms \cite{AW70,B78}, where entropy is again
a fundamental invariant. In mathematical physics SFTs are often called
hard-core models, and are used to model a wide variety of physical
systems; this is the thermodynamic formalism \cite{R04}. In this
setup it is of central importance to understand the equilibrium states
of the system, which in are the invariant measures of maximal entropy.

It is well known that in the one-dimensional case the entropy of an
SFT may be effectively calculated, since it is the logarithm of the
spectral radius of a certain positive integer matrix which is derived
from the combinatorial description of the system. D. Lind \cite{lind84}
has given an algebraic characterization of the numbers which arise
as entropies of one-dimensional SFTs. A Perron number is a real algebraic
integer greater than $1$ and greater than the modulus of its algebraic
conjugates. The entropies of one-dimensional SFTs are precisely the
non-negative rational multiples of logarithms of Perron numbers.

In higher dimensions the problem becomes much more difficult. The
dynamics of multidimensional SFTs is vastly more complicated than
their one-dimensional counterparts. For instance, irreducible multidimensional
SFTs may have more than one measure of maximal entropy \cite{burton94},
and zero entropy can coexist with rather complex dynamics \cite{mozes89}.
In general it is undecidable whether a given set of rules define a
nonempty SFT \cite{B66,robinson71}. Regarding the entropy, even when
the rules defining an SFT enjoy good symmetry properties, calculating
the entropy is usually beyond current technology. As a result numerical
methods have been developed to approximate the entropy (e.g. \cite{F05}),
but these usually apply to  a restricted class systems. 

One should note that for certain $\mathbb{Z}^{d}$-actions which arise
as automorphisms of compact groups (a class which includes some SFTs),
explicit expressions for the entropy have been obtained by D. Lind,
K. Schmidt and T. Ward \cite{LS90}. We note however that while these
expressions are explicit they do not provide much information on the
properties of the entropies, e.g. whether they are algebraic, well
approximable, etc.

In this paper we characterize those real numbers which can occur as
entropies of multidimensional SFTs in terms of their computation-theoretic
properties. It is natural to say that a real number $h$ is \emph{computable}
if it can be calculated to any desired accuracy. More precisely, $h$
is \emph{}computable if there is an algorithm which, given input $n\in \mathbb{N}$,
produces a rational number $r(n)$ with $|h-r(n)|<1/n$. For example,
every algebraic number is computable (since there are numerical methods
for computing the roots of an integer polynomial), and so are $e,\pi $,
since they can be written as power series with computable coefficients
and rate of convergence. 

A weaker notion is the following. A real number $h$ is \emph{right
recursively enumerable} (sometimes called \emph{upper semi recursive)}
if there exists a Turing machine which, given $n$, computes a rational
number $r(n)\geq h$ such that $r(n)\rightarrow h$ (equivalently,
the right Dedekind cut $\{q\in \mathbb{Q}\, :\, q>h\}$ is a recursive
set of rationals). 

The class of right recursively enumerable numbers is countable since
algorithms may be put in one-to-one correspondence with finite $0,1$-valued
sequences, and hence there are only countably many of them. If $h$
is computable then there is an algorithm computing $r(n)$ with $|h-r(n)|<\frac{1}{n}$,
so the computable sequence $r(n)+\frac{1}{n}$ converges to $h$ from
above. This shows that the class of right recursively enumerable numbers
contains the computable numbers, and it can be shown to be strictly
larger. For more information, see \cite{KO91}.

\begin{thm}
\label{thm:main}For $d\geq 2$ the class of entropies of $d$-dimensional
SFTs is the class of non-negative right recursively enumerable numbers. 
\end{thm}
The property of right recursive enumerability is a necessary condition
for a number to be the entropy of an SFT because the naive approximation
algorithm, which counts locally admissible patterns on cubes, converges
from above to the entropy. This follows from the work of Friedland
\cite{friedland97}; we provide a different proof below. The main
novelty here is the sufficiency of the condition. 

A sofic system is a factor of an SFT, i.e. an ensemble of colorings
of $\mathbb{Z}^{d}$ obtained from a fixed SFT $X$ by applying a
local transformation to each coloring in $X$ (for a definition, see
section \ref{sec:Preliminaries}). In the one dimensional case, Coven
and Paul \cite{CP75} showed that every sofic system $Y$ can be extended
to an SFT $X$ with of the same entropy as $Y$. In particular, this
implies that the class of entropies of sofic shifts is the same as
that of SFTs. Whether the covering theorem is true in the multidimensional
case is still open and seems quite hard (see \cite{D06} for a partial
result). However some circumstantial evidence in favor of the covering
theorem is provided by the following:

\begin{thm}
\label{thm:sofic}For $d\geq 2$, the class of entropies of $d$-dimensional
sofic shifts is the same as that of $d$-dimensional SFTs. 
\end{thm}
This is a consequence of the fact that the entropy of sofic shifts
is right recursively enumerable (corollary \ref{cor:sofic-entropies}
below), and the fact that an SFT is in particular a sofic system.

It is worth emphasizing that since there are non-computable numbers
which are right recursively enumerable, it follows from theorem \ref{thm:main}
that there are SFTs whose entropy cannot be computed effectively (a
similar situation is known for cellular automata \cite{Hurd92} and
general subshifts \cite{S06}). However, if one assumes strong enough
mixing properties of the system the situation improves. Recall that
an SFT is \emph{irreducible} if any two admissible patterns far enough
apart may be extended to the whole lattice (see section \ref{sec:Preliminaries}). 

\begin{thm}
\label{thm:irreducible}The entropy of an irreducible SFT is computable.
\end{thm}
We do not know if this condition is also sufficient. 

The rest of this paper is organized as follows. In the next section
we introduce some notation and background. In section \ref{sec:Entropies-of-SFTs-are-co-RE}
we prove that the entropy of any SFT or sofic shift is right recursively
enumerable, and that of an irreducible SFT is computable. In section
\ref{sec:Outline} we outline the construction which constitutes the
proof of the other direction of theorem \ref{thm:main}. Sections
\ref{sec:base}--\ref{sec:calculating-entropy} give the details of
the construction. In section \ref{sec:open-problems} we discuss discuss
some open problems.

\begin{acknowledgement*}
This work was done during the authors' graduate studies, and we would
like to thank our advisors, Benjamin Weiss and Jon Aaronson, for their
support and advice. We also thank Mike Boyle for his comments.
\end{acknowledgement*}

\section{\label{sec:Preliminaries}Preliminaries}

In this section we provide some background from symbolic dynamics
and define SFTs and entropy. See \cite{LM95,K98} for more information
on these subjects.

\subsection{Subshifts and SFTs.}

For a finite alphabet $\Sigma $ let $\Sigma ^{\mathbb{Z}^{d}}$ be
the space of $\Sigma $-colorings of $\mathbb{Z}^{d}$ (this is called
the full shift on $\Sigma $). For a subset $F\subseteq \mathbb{Z}^{d}$
we refer to a function $a\in \Sigma ^{F}$ as a coloring of $F$ or
an $F$-pattern. We say that patterns $a\in \Sigma ^{F}$ and $b\in \Sigma ^{F+u}$
are \emph{congruent} if $a(v)=b(v+u)$ for every $v\in F$. We say
that a pattern $a\in \Sigma ^{F}$ appears at $u$ in a pattern $b\in \Sigma ^{E}$
if $b|_{F+u}$ and $a$ are congruent.

If $E\subseteq F$ and $a\in \Sigma ^{F}$ then $a$ induces a coloring
of $E$ by restriction, namely $a|_{E}$. For a finite set $F\subseteq \mathbb{Z}^{d}$
and pattern $a\in \Sigma ^{F}$ the cylinder set defined by $a$ is\[
[a]=\{x\in \Sigma ^{\mathbb{Z}^{d}}\, :\, x|_{F}=a\}\]
 We endow $\Sigma ^{\mathbb{Z}^{d}}$ with the product topology, which
is generated by the cylinder sets and makes $\Sigma ^{\mathbb{Z}^{d}}$
into a compact metrizable space.

For $u\in \mathbb{Z}^{d}$ let $\sigma ^{u}:\Sigma ^{\mathbb{Z}^{d}}\rightarrow \Sigma ^{\mathbb{Z}^{d}}$
be the homeomorphisms \[
(\sigma ^{u}(x))(v)=x(v+u)\qquad v\in \mathbb{Z}^{d}\]
 This gives an action of $\mathbb{Z}^{d}$ on $\Sigma ^{\mathbb{Z}^{d}}$
called the shift action. A subset $X\subseteq \Sigma ^{\mathbb{Z}^{d}}$
is invariant under the shift action if $\sigma ^{u}(X)=X$ for every
$u\in \mathbb{Z}^{d}$. A closed invariant set $X\subseteq \Sigma ^{\mathbb{Z}^{d}}$
is called a $\mathbb{Z}^{d}$-subshift.

A $d$-dimensional subshift of finite type (SFT) is defined by a finite
alphabet $\Sigma $, a finite set $F\subseteq \mathbb{Z}^{d}$, and
a collection $L\subseteq \Sigma ^{F}$ of $\Sigma $-colorings of
$F$, called the \emph{syntax}. A $\Sigma $-coloring $x\in \Sigma ^{\mathbb{Z}^{d}}$
of $\mathbb{Z}^{d}$ is \emph{admissible} for $L$ if the pattern
induced by $x$ on every translate of $F$ is congruent to a pattern
in $L$. The SFT defined by $L$ is the set $X\subseteq \Sigma ^{\mathbb{Z}^{d}}$
of all admissible $x$. From the definition it is clear that an SFT
is closed and shift-invariant. 

Given an SFT $X$ defined by a syntax $L\subseteq \Sigma ^{F}$, we
say that a finite pattern is \emph{globally admissible} for $X$ if
it appears in $X$. In contrast we say that a pattern $a\in \Sigma ^{E}$
is \emph{locally admissible} if $a|_{F+u}$ is congruent to a pattern
in $L$ whenever $F+u\subseteq E$. A globally admissible pattern
is locally admissible, but the latter is not true in general.

An SFT $X\subseteq \Sigma ^{\mathbb{Z}^{d}}$ is \emph{irreducible}
if there is a constant $r>0$, called a \emph{gap}, such that for
every $A,B\subseteq \mathbb{Z}^{d}$ satisfying $\left\Vert u-v\right\Vert _{\infty }\geq r$
for $u\in A\, ,\, v\in B$, and for every pair of globally admissible
$a\in \Sigma ^{A}$ and $b\in \Sigma ^{B}$, there is a point $x\in X$
with $x|_{E}=a$ and $x|_{B}=b$ (in other words, $a\cup b$ is globally
admissible).

\subsection{Topological entropy of subshifts.}

For a subshift $X\subseteq \Sigma ^{Z^{d}}$ and $F\subseteq \mathbb{Z}^{d}$
we say that a pattern $a\in \Sigma ^{F}$ appears in $X$ if $a=x|_{F}$
for some $x\in X$. For a set $F$ let $N_{X}(F)$ denote the number
of distinct $\Sigma $-colorings of $F$ which appear in $X$. Let
\[
F_{n}=\{1,\ldots ,n\}^{d}\]
 denote the discrete $d$-dimensional cube of side $n$. The (topological)
entropy $h(X)$ of $X$ is defined by\[
h(X)=\lim _{n\rightarrow \infty }\frac{1}{|F_{n}|}\log N_{X}(F_{n}).\]
 By convention the logarithm is to base $2$. The limit above exists,
and is in fact equal to $\inf _{n\in \mathbb{N}}\frac{1}{|F_{n}|}\log N_{X}(F_{n})$.

\subsection{Products, factors and isomorphism}

Let $X\subseteq \Sigma ^{\mathbb{Z}^{d}}$ and $Y\subseteq \Delta ^{\mathbb{Z}^{d}}$
be two $\mathbb{Z}^{d}$-subshifts. The product system $X\times Y\subseteq (\Sigma \times \Delta )^{\mathbb{Z}^{d}}$
is then a symbolic $\mathbb{Z}^{d}$-system also, and satisfies $h(X\times Y)=h(X)+h(Y)$.

A continuous onto map $\varphi :X\rightarrow Y$ is called a factor
map if it commutes with the action, i.e. $\sigma ^{u}\circ \varphi =\varphi \circ \sigma ^{u}$
for all $u\in \mathbb{Z}^{d}$. An isomorphism is an invertible factor
map. Both entropy and the property of being an SFT are invariants
of isomorphism (although isomorphic SFTs are usually not defined by
the same syntax), as is irreducibility.

A factor of an SFT is called a sofic system. In general a sofic system
is not an SFT .

Every factor map $\varphi :X\rightarrow Y$ arises from a so-called
block code, which means the following: There exists a finite set $F\subseteq \mathbb{Z}^{d}$
and a function $\varphi _{0}:\Sigma ^{F}\rightarrow \Delta $ such
that\[
(\varphi (x))(u)=\varphi _{0}((\sigma ^{u}x)|_{F})\]
 Conversely, given such a $\varphi _{0}$ we can define $\varphi $
by this formula, and then $\varphi $ is a factor map from $X$ onto
its image.

A factor map $\pi :X\rightarrow Y$ of symbolic systems $X\subseteq \Sigma ^{\mathbb{Z}^{d}}$
and $Y\subseteq \Delta ^{\mathbb{Z}^{d}}$ is called a one-block map
if it is determined by a single symbol, i.e. it is induced by a map
$\varphi _{0}:\Sigma \rightarrow \Delta $. We will always assume
our factor maps are one-block maps. There is no loss of generality
in this since given a factor map $\varphi :X\rightarrow Y$ there
is a system $X'$ isomorphic to $X$ via an isomorphism $\psi :X'\rightarrow X$
so that the factor map $\psi \varphi :X'\rightarrow Y$ is a one-block
map. 

Similarly, an SFT is called one-step if it is defined by a syntax
$L\subseteq \Sigma ^{\{0,1\}^{d}}$. Every SFT is isomorphic to a
one-step SFT. Note that for a one-step SFT, if $a\in \Sigma ^{F_{n}}$
and $b\in \Sigma ^{\mathbb{Z}^{d}\setminus F_{n-1}}$ are globally
admissible patterns and they agree on the boundary of $F_{n}$ (i.e.
$a|_{F_{n}\setminus F_{n-1}}=b|_{F_{n}\setminus F_{n-1}}$) then $a\cup b$
is globally admissible.

\subsection{Invariant measures and entropy. }

Given a symbolic system $X$, a Borel measure $\mu $ on $X$ is invariant
under the shift action if $\mu (\sigma ^{u}(A))=\mu (A)$ for every
Borel set $A\subseteq X$ and every $u\in \mathbb{Z}^{d}$. We denote
the set of invariant Borel probability measures by $\mathcal{M}(X)$.
The weak-{*} topology on $\mathcal{M}(X)$ is the topology in which
$\mu _{n}\rightarrow \mu $ if $\int fd\mu _{n}\rightarrow \int fd\mu $
for every continuous function $f$ on $X$. This makes $\mathcal{M}(X)$
into a compact metrizable space. 

For $\mu \in \mathcal{M}(X)$ we denote by its measure-theoretic entropy
by $h(\mu )$. We recall the following facts:

\begin{enumerate}
\item The entropy function $h:\mathcal{M}(X)\rightarrow \mathbb{R}^{+}$
is upper semi-continuous.
\item The variational principle: $h(X)=\max _{\mu \in M(X)}h(\mu )$. 
\end{enumerate}
See \cite{MCK76} for definitions, proofs and a detailed discussion
of the one-dimensional case, or \cite{M76} for a proof of the variational
principle in the multidimensional case.

\section{\label{sec:Entropies-of-SFTs-are-co-RE}Computability of Entropies}

In this section we show that the entropy of an SFT is right recursively
enumerable. This follows from the work of Friedland \cite{friedland97},
but for completeness we give an short alternative proof and extend
the result to sofic systems. We also prove that the entropy of an
irreducible SFT is computable.

Let the syntax $L\subseteq \Sigma ^{F}$ define a (possibly empty)
SFT $X$. The definition of entropy provides us with the sequence
$N_{n}=N_{X}(F_{n})$ such that $\frac{1}{n^{d}}\log N_{n}$ which
converges to $h(X)$ from above, and if $N_{n}$ is computable this
sequence shows that $h$ is upper-semi recursive. However, $N_{n}$
is not computable in general. Indeed, determining whether $N_{n}>0$
is equivalent to deciding if the SFT defined by $L$ is nonempty,
and this is in general undecidable \cite{robinson71,B66}. 

Let us say a finite pattern $a\in \Sigma ^{F_{n}}$ is locally admissible
if $a|_{F+u}$ is congruent to a pattern in $L$ whenever $F+u\subseteq F_{n}$.
Instead of $N_{n}$, consider the sequence\[
\tilde{N}_{n}=\#\{\textrm{locally admissible }F_{n}\textrm{-patterns}\}\]
Clearly $\tilde{N}_{n}$ is computable. If $x\in X$ then $x|_{F_{n}}$
is one of the patterns counted by $\tilde{N}_{n}$, so $\tilde{N}_{n}\geq N_{n}$.
The inequality can be strict, because not all locally admissible $F_{n}$-patterns
need arise in this way: there can be locally admissible finite patterns
which don't extend to globally admissible coloring of $\mathbb{Z}^{d}$.
Nonetheless,

\begin{thm}
\label{thm:SFT-entropies-are-co-RE}For $L,X$ and $\tilde{N}_{n}$
as above, $\frac{1}{n^{d}}\log \tilde{N}_{n}\rightarrow h(X)$ from
above. Consequently, $h(X)$ is right recursively enumerable. 
\end{thm}
\begin{proof}
Denote \[
\tilde{h}=\limsup _{n\rightarrow \infty }\frac{1}{n^{d}}\log \tilde{N}_{n}\]
 Since $\widetilde{N}_{n}\geq N_{n}$ and $\frac{1}{n^{d}}\log N_{n}\geq h(X)$,
we have $\frac{1}{n^{d}}\log \widetilde{N}_{n}\geq h(X)$, so it suffices
to show that $\tilde{h}\leq h(X)$.

Define a sequence of measures $\nu _{n}$ on $\Sigma ^{\mathbb{Z}^{d}}$
as follows. Let $W_{n}\subseteq \Sigma ^{F_{n}}$ be the set of locally
admissible colorings of $F_{n}$. Let $\mu _{n}$ denote the probability
measure obtained by coloring each translate $F_{n}+u$ for $u\in n\mathbb{Z}^{d}$
independently and uniformly with patterns from $W_{n}$. Let $\nu _{n}=\sum _{u\in F_{n}}\sigma ^{u}\mu _{n}$.
Then $\nu _{n}$ is an invariant probability measure and its entropy
is easily shown to be \[
h(\nu _{n})=\frac{1}{n^{d}}\log \tilde{N}_{n}\]
 Let $\nu _{n(k)}$ be a subsequence such that $h(\nu _{n(k)})\rightarrow \tilde{h}$
and let $\nu $ be a weak-{*} accumulation point $\nu _{n(k)}$; we
may assume $\nu _{n(k)}\rightarrow \nu $. Since entropy is upper
semi-continuous in the weak-{*} topology, we have \[
h(\nu )\geq \tilde{h}\]

On the other hand we claim that $\nu (X)=1$, so $\nu $ can be regarded
as an invariant probability measure on $X$. To show this, we prove
that $\nu ([a])=0$ for any $a\in \Sigma ^{F}\setminus L$ , where
$[a]$ is the cylinder set defined by $a$. Indeed, for every $k$
and $u\in F_{k}$, if $(F+u)\subset F_{k}$ then $\mu _{k}(\sigma ^{-u}([a]))=0$,
so \[
\mu _{k}([a])\leq \frac{1}{k^{d}}\#\{u\in F_{k}:\; (F+u)\not \not \subseteq F_{k}\}\leq \frac{k^{d}-(k-\diam F)^{d}}{k^{d}}\]
 where $\diam F$ is the diameter of $F$ with respect to the norm
$\left\Vert u\right\Vert _{\infty }=\max _{i=1\ldots d}|u_{i}|$.
It is easy to see that $\nu _{k}([a])=\mu _{k}([a])$, so\[
\nu ([a])=\lim _{k\rightarrow \infty }\nu _{n(k)}([a])=0\]
 Finally, the variational principle implies that $h(\nu )\leq h(X)$,
and the theorem follows. 

\end{proof}
With the same notation as above, let $Y\subseteq \Delta ^{\mathbb{Z}^{d}}$
be a symbolic factor of $X$ arising from a one-block map $\varphi _{0}:\Sigma \rightarrow \Delta $
and its pointwise extension $\varphi :X\rightarrow Y$. Write \[
\tilde{M}_{n}=\#\{\varphi (a)\in \Delta ^{F_{n}}:\; a\in W_{n}\}\]
 where as before, $W_{n}$ is the set of locally admissible $F_{n}$-patterns
for $L$.

\begin{thm}
\label{thm:sofic-entropies-are-co-RE}With the above notation, $\frac{1}{|F_{n}|}\log \tilde{M}_{n}\rightarrow h(Y)$
from above. Consequently $h(Y)$ is right recursively enumerable. 
\end{thm}
\begin{proof}
Denote $\tilde{h}(Y)=\limsup \frac{1}{n^{d}}\log \widetilde{M}_{n}$.
Since $\varphi $ is onto we have $\tilde{M}_{n}\geq N_{Y}(F_{n})$,
so $\frac{1}{n^{d}}\log \widetilde{M}_{n}\geq h(Y)$. Thus we only
need to show $\tilde{h}(Y)\leq h(Y)$.

Let $\theta _{k}$ be measures on $\Delta ^{\mathbb{Z}^{d}}$ defined
by coloring each translate $F_{k}+u$ for $u\in n\mathbb{Z}^{d}$
with patterns drawn uniformly from $\{\varphi (a)\, :\, a\in W_{k}\}$.
Then $\eta _{k}=\frac{1}{k^{2}}\sum _{u\in F_{k}}\sigma ^{u}\theta _{k}$
is an invariant measure on $\Delta ^{\mathbb{Z}^{d}}$, and $h(\eta _{k})=\frac{1}{k^{d}}\log \tilde{M}_{k}$.
Let $\mu _{k}$ be measures on $\Sigma ^{\mathbb{Z}^{d}}$ such that
the pattern on $F_{k}+u$ for $u\in k\mathbb{Z}^{d}$ is drawn from
$W_{k}$ according to a distribution which projects under $\varphi $
to the uniform distribution on $\{\varphi (a)\, :\, a\in W_{k}\}$.
Thus $\theta _{k}=\varphi (\mu _{k})$. Let $\nu _{k}=\sum _{u\in F_{k}}\sigma ^{u}\mu _{k}$,
so that $\eta _{k}=\varphi (\nu _{k})$. Choose a subsequence $\eta _{n(k)}$
so that there is a measure $\eta $ on $Y$ with $\eta _{n(k)}\rightarrow \eta $
and there is a measure $\nu $ on $X$ with $\nu _{n(k)}\rightarrow \nu $;
so $\eta =\varphi (\nu )$ satisfies $\eta _{n(k)}\rightarrow \eta $.
By upper semi-continuity, $h(\eta )\geq \lim h(\eta _{n(k)})=\tilde{h}$,
and so we will be done if we show that $\eta $ is supported on $Y$.
For this it is enough to show that $\nu $ is supported on $X$, i.e.
that $\nu ([b])=0$ whenever $b\in \Sigma ^{F}\setminus L$. The proof
of this is identical to the proof of the same statement at the end
of theorem \ref{thm:SFT-entropies-are-co-RE}. 
\end{proof}
\begin{cor}
\label{cor:sofic-entropies}The entropy of every sofic shift is right
recursively enumerable. 
\end{cor}
\begin{proof}
As noted in section \ref{sec:Preliminaries}, every sofic shift is
a one-block factor of some SFT. 
\end{proof}
We turn now to irreducible SFTs and the proof of theorem \ref{thm:irreducible}.
Let $Q_{n}=\{-n,-n+1,\ldots ,n-1,n\}^{d}$ denote the symmetric cube;
note that $N_{X}(Q_{n})=N_{X}(F_{2n+1})$. For an SFT $X$ we say
that globally admissible patterns $a\in \Sigma ^{Q_{k}}$ and $b\in \Sigma ^{Q_{n}}$
are $r$-compatible if $n\geq k+r+1$ and $a\cup (b|_{Q_{n}\setminus Q_{k+r}})$
is globally admissible. Note that if $r<s$ and $a\in \Sigma ^{Q_{k}}$
and $b\in \Sigma ^{Q_{n}}$ are $r$-compatible, then they are $s$-compatibe,
provided $n\geq k+s+1$. Clearly if $X$ is irreducible with gap $r$
then every two such patterns are $r$-compatible. 

\begin{lem}
Let $X\subseteq \Sigma ^{\mathbb{Z}^{d}}$ be a non-empty irreducible
SFT and $a\in \Sigma ^{Q_{k}}$. Then precisely one of the following
holds:
\begin{enumerate}
\item $a\neq b|_{Q_{k}}$ for every large enough $N$ and every locally
admissible $b\in \Sigma ^{Q_{N}}$.
\item For every large enough $N$ and locally admissible $b\in \Sigma ^{Q_{N}}$,
the patterns $a,b$ are $\sqrt{N}$-compatible.
\end{enumerate}
In the first case $a$ is not globally admissible; in the second it
is.
\end{lem}
\begin{proof}
Fix $a\in \Sigma ^{Q_{k}}$. By compactness, $a$ is not globally
admissible if and only if (1) holds.

Suppose now that $a$ is globally admissible for $X$ and let $r$
be a gap for $X$. We may assume that $X$ is a one-step SFT. For
every $b\in \Sigma ^{Q_{k+r+1}}$, if $a,b$ are \emph{not} $r$-compatible
then by irreducibility $b$ is not globally admissible. Hence by part
(1) we see that for large enough $N$ the pattern $b$ does not appear
at the origin in any locally admissible $c\in \Sigma ^{Q_{N}}$. Since
there are finitely many such $b$'s, we see that $a$ and $c$ are
$r$-compatible for large enough $N$ and for locally admissible $c\in \Sigma ^{F_{N}}$.
Since $\sqrt{N}>r$ eventually, this implies the (2). 

Finally, (2) implies that $a$ is admissible by irreducibility and
tha fact that $X$ is non-empty.
\end{proof}
\begin{cor}
For a non-empty irreducible SFT $X$ it is decidable whether a finite
pattern $a$ is globally admissible.
\end{cor}
\begin{proof}
To decide if $a$ is globally admissible, find the first $N$ for
which one of the conditions of the proposition holds (the conditions
are finitely checkable). If (2) holds then $a$ is globally admissible;
otherwise it is not. Note that to apply the proposition one does not
need to know the gap.
\end{proof}
A number $h$ is left recursively enumerable if there is an algorithm
which, given $n$, produces a rational number $s(n)$ with $s(n)\rightarrow h$
and $s(n)\leq h$. If $h$ is both right and left recursively enumerable
then it is computable. To see this let $r(n),s(n)$ be computable
sequences with $s(n)\leq h\leq r(n)$ and $r(n),s(n)\rightarrow h$.
Now given $n$, we can calculate $r(k),s(k)$ for $k=1,2,3\ldots $
until such a $k$ is reached that $r(k)-s(k)<\frac{1}{n}$. Then $r(k)$
satisfies $|r(k)-h|<\frac{1}{n}$. This algorithm shows that $h$
is computable. 

We can now prove theorem \ref{thm:irreducible}, which we repeat here
for convenience:

\begin{thm*}
The entropy of an irreducible SFT is computable.
\end{thm*}
\begin{proof}
Let $X$ be an irreducible SFT, and we may assume it is non-empty.
We already know that $h(X)$ is right recursively enumerable, so it
suffices to show that it is left recursively enumerable, i.e. to exhibit
an algorithm which given $n\in \mathbb{N}$ returns a rational number
$s(n)$ such that $s(n)\rightarrow h(X)$ and $s(n)\leq h$. 

The algorithm is as follows. First, identify all the globally admissible
patterns $a_{1},\ldots ,a_{k(n)}\in \Sigma ^{Q_{n}}$ (this is computable
by corollary above). With this notation we have $k(n)=N_{X}(Q_{n})$
and $\frac{1}{|Q_{n}|}\log k(n)\rightarrow h(X)$. Next, find the
smallest number $r'$ so that each globally admissible pattern $b\in \Sigma ^{Q_{n+r'+1}}$
is $r'$-compatible with $a_{i}$ for $i=1,\ldots ,k(n)$. Set\[
s(n)=\frac{1}{|Q_{n+r'}|}\log k(n)\]

Note that $r'\leq r$, where $r$ is a fixed gap for $X$. Hence \[
s(n)\geq \frac{|Q_{n}|}{|Q_{n+r}|}\cdot \frac{1}{|Q_{n}|}\log k(n)\rightarrow h(X)\]
On the other hand, consider a large $Q_{m}$, and consider the collection
of translates of $Q_{n}$ by elements of the lattice $2(n+r')\mathbb{Z}^{d}$
which fall inside $Q_{m}$. By choice of $r'$ we can color each of
these translates in an arbitrary globally admissible way and complete
it to a globally admissible $Q_{m}$ pattern. Since the number of
translates is $\frac{1}{|Q_{n+r'}|}|Q_{m}|$ (for convenience assume
that $m$ is a multiple of $2(n+r'$), we see that \[
k(m)\geq k(n)^{\frac{1}{|Q_{n+r'}|}\cdot |Q_{m}|}\]
so letting $m\rightarrow \infty $,\[
s(n)=\frac{1}{|Q_{n+r'}|}k(n)\leq \frac{1}{|Q_{m}|}\log k(m)\rightarrow h(X)\]
 hence $s(n)\leq h(X)$, and also $s(n)\rightarrow h(X)$, as desired.
\end{proof}
Note that the algorithm given in the proof does not require prior
knowledge of a gap for the $X$. It may of course be applied to any
SFT, but in that case may not halt on some inputs, and even if it
does the sequence $s(n)$ will not necessarily behave as above.

\section{\label{sec:Outline}Outline of the main construction}

Let $h$ be a right recursively enumerable number. To prove the remaining
direction of theorem \ref{thm:main}, we must construct for every
$d\geq 2$ a $d$-dimensional SFT with entropy $h$. We first make
some simplifying assumptions. We may restrict ourselves to dimension
$2$, since given an SFT $X\subseteq \Sigma ^{\mathbb{Z}^{d}}$ the
system $X'\subseteq \Sigma ^{\mathbb{Z}^{d+1}}$ defined by\[
X'=\{x'\in \Sigma ^{\mathbb{Z}^{d+1}}\, :\, \forall j\in \mathbb{Z}\, \exists x\in X\, \forall u\in \mathbb{Z}^{d}\, x'(u,j)=x(u)\}\]
 is easily seen to be a $d+1$-dimensional SFT and $h(X')=h(X)$.
Furthermore, since (a) the product of SFTs is an SFT, (b) $h(X\times Y)=h(X)+h(Y)$
and (c) $n$ is the entropy of the full shift on $2^{n}$ symbols,
it suffices to prove the statement under the assumption that $h\in [0,1]$.

Our construction has three main steps:

\begin{description}
\item [Step~1~:~Constructing~the~base~(section~\ref{sec:base})]We
construct an SFT $X$ some of whose symbols are marked $0,1$, and
such that the density of $1$'s in each point of $X$ is very uniform.
It will be possible to estimate this density by observing any sufficiently
large and well-distributed set of coordinates.
\item [Step~2~:~Pruning~(section~\ref{sec:Pruning})]In this step
we {}``kill'' all points $x\in X$ such that the frequency of $1$'s
in $x$ is strictly greater then $h$. In this way we obtain an SFT
$Y$ such that the symbol $1$ appears in each $y\in Y$ with frequency
at most $h$, and for some points the frequency is $h$. Furthermore,
$Y$ will still have zero entropy. We achieve this by superimposing
another layer on top of $X$ which represents calculations of a certain
Turing machine, using as input the underlying patterns from $X$.
This machine halts when it detects a density of $1$'s greater than
$h$. The result is that a point $x\in X$ with density of $1$'s
greater than $h$ cannot be extended to a pattern in $Y$; otherwise,
it can be. 
\item [Step~3~:~Adding~{}``Random''~bits~(section~\ref{sec:calculating-entropy})]We
extend $Y$ to an SFT $Z$ by allowing two new symbols, say ''$\alpha $''
and {}``$\beta $'', to appear independently over every occurrence
of a $1$ in $Y$. This system $Z$ has entropy $h$.
\end{description}
For steps 1 and 2 we utilize certain SFTs with special geometric and
arithmetic properties. The existence of such systems, and their use
in representing Turing machines in SFTs, appears first in Robinson's
paper \cite{robinson71}. However, we will not refer directly to Robinson's
construction, which would in any case require some modification to
suit our needs. Instead we rely on a theorem of Mozes \cite{mozes89}
about the realization of substitution systems by SFTs. This theorem,
which allows us to easily construct variants of Robinson's system,
is presented in the next section together with another technical definition.
Following that we give the details of steps 1, 2 and 3.

Before moving on, we note that our arguments give the following result,
which may be of independent interest:

\begin{thm}
A real number $r\geq 0$ is right recursively enumerable if and only
if there is an alphabet $\Sigma $, a symbol $a\in \Sigma $ and an
SFT $X\subseteq \Sigma ^{\mathbb{Z}^{d}}$ such that \[
\sup _{x\in X}\, \lim _{n\rightarrow \infty }\frac{1}{|F_{n}|}\#\{u\in F_{n}\, :\, x(u)=a\}=r\]
(and in particular the limit above exists for every $x\in X$). Furthermore
if $r$ is computable then one can find $\Sigma ,a,X$ so that $\lim _{n\rightarrow \infty }\frac{1}{|F_{n}|}\#\{u\in F_{n}\, :\, x(u)=a\}=r$
for every $x\in X$. 
\end{thm}

\section{Substitutions and superpositions}

In this section we describe two technical devices for constructing
SFTs.

\subsection{\label{sub:Substitution-systems}Subshifts defined by Substitution }

Given a finite alphabet $\Sigma $, a \emph{substitution rule} is
a map $s:\Sigma \rightarrow \Sigma ^{F_{k}}$ for some integer $k>1$,
where $F_{k}=\{1,\ldots ,k\}\times \{1,\ldots ,k\}$ (in the terminology
of \cite{mozes89}, this is a deterministic $k\times k$ substitution
system with property $A$). The map $s$ extends naturally to a map
$s_{n}:\Sigma ^{F_{n}}\rightarrow \Sigma ^{F_{n\cdot k}}$ by identifying
$\Sigma ^{F_{n\cdot k}}$ with $(\Sigma ^{F_{k}})^{F_{n}}$. 

Starting from a single symbol located at $(1,1)\in \mathbb{Z}^{2}$
and iterating the substitution map, we obtain a sequence of colorings
of $F_{k^{n}}$ for $n=0,1,2\ldots $. Such patterns are called $s$-blocks.
A point $x\in \Sigma ^{\mathbb{Z}^{2}}$ is admissible for $s$ if
every finite subpattern of $x$ appears in some $s$-block. The subshift
$W\subseteq \Sigma ^{\mathbb{Z}^{2}}$ associated with $s$ is the
set of admissible patterns; this is seen to be closed and shift invariant.

Define $s_{\infty }:W\rightarrow W$ by applying $s$ to each symbol
of $x$; more precisely, $s_{\infty }(x)(u)=s(x(u'))(u'')$, where
$u'\in \mathbb{Z}^{2}$ and $u''\in F_{k}$ are the unique vectors
such that $u=ku'+u''$. Clearly $s_{\infty }$ maps $W$ into $W$.
We say that $x$ is derived from $y$ if $T^{v}x=s_{\infty }(y)$
for some $v\in F_{k}$. It is not hard to show that each $x\in W$
is derived from some $y\in W$; if this $y$ is unique, we say that
$s$ has \emph{unique derivation}. 

\begin{thm}
\label{thm:subst}(Theorem 4.5 of \cite{mozes89}) Let $s:\Sigma \rightarrow \Sigma ^{F_{k}}$
be a substitution rule with unique derivation and let $W$ be the
associated dynamical system. Then there exists an alphabet $\Delta $,
an SFT $\widetilde{W}\subseteq \Delta ^{\mathbb{Z}^{2}}$, and a one-block
factor map $\varphi :\widetilde{W}\rightarrow W$. Furthermore $\varphi $
is an injection on a set having full measure with respect to every
invariant measure on $\widetilde{W}$.
\end{thm}
Note that theorem \ref{thm:subst} is false in dimension $d=1$. 

\begin{prop}
If $s,W$ and $\widetilde{W}$ are as in the theorem then $h(\widetilde{W})=0$.
\end{prop}
\begin{proof}
For any $\mu $ invariant on $\widetilde{W}$, the map $\varphi $
is an isomorphism of dynamical systems between $(\widetilde{W},\mu )$
and $(W,\varphi \mu )$ where $\varphi \mu $ is the push-forward
of $\mu $ to $\widetilde{W}$. Hence it suffices to show that the
latter system has zero measure-theoretic entropy. By the variational
principle it suffices to show that $h(W)=0$. Fix $m$. Since every
large enough $s$-blocks is composed of an array of smaller $s$-blocks
of dimension $k^{m}\times k^{m}$ arranged in a square, it follows
that for $n>k^{m}$ an admissible $F_{n}$-pattern can be decompose
$F_{n}$ into $([\frac{n}{k^{m}}]-2)^{2}$ disjoint $s$-blocks of
dimension $k^{m}\times k^{m}$ together with a {}``small'' remaining
region near the boundary. Thus the number of $F_{n}$ patterns is
at most \[
N_{F_{n}}(W)\leq \#\{k\times k\textrm{ }s\textrm{-blocks}\}^{([n/k^{m}]-2)^{2}}\cdot |\Sigma |^{4nk^{m}}\]
where the second term on the right hand side is the number of ways
to fill in the region near the boundary of $F_{n}$ not covered by
the $s$-blocks. Since there are only $|\Sigma |$ different $s$-blocks
of dimension $k^{m}\times k^{m}$ (because each is derived from one
of the original symbols), for all large enough $n$ we have\[
\frac{1}{n^{2}}\log N_{F_{n}}(W)\leq \frac{([\frac{n}{k^{m}}]-2)^{2}\log |\Sigma |}{n^{2}}+\frac{4k^{m}\log |\Sigma |}{n}\rightarrow \frac{1}{k^{m}}\]
as $m$ was arbitrary, $h(W)=0$.
\end{proof}
We use theorem \ref{thm:subst}, which is due to Mozes, to construct
systems similar in many respects to Robinson's system from \cite{robinson71}.
We remark that although this allows a more economical exposition the
gain is cosmetic. Indeed, the proof of theorem \ref{thm:subst} relies
on an elaborate extension of Robinson's techniques. There has recently
been a revival of interest in substitutions and their realization
using local rules; see e.g. \cite{GS98}.

\subsection{Superposition}

Given an SFT $X$ defined by a syntax $L$, \emph{superposition} is
a syntactic process which gives an SFT $X'$ which factors onto a
subshift of $X$. Informally, this is done by adding data to each
symbol of $X$ and enriching the syntax with rules relating to this
new data. 

More precisely, suppose $X$ is an SFT defined by a syntax $L\subseteq \Sigma ^{F}$.
A system $Y$ is superimposed over $X$ if it is obtained by the following
process. (a) For a finite set $\Delta $, we replace each symbol of
$\sigma \in \Sigma $ with one or more symbols if the form $(\sigma ,\delta )\in \Sigma \times \Delta $.
Let $\Sigma '$ be the set of these pairs. For the new symbol $(\sigma ,\delta )\in \Sigma '$,
we say that $\delta $ is superimposed over $\sigma $; we also frequently
refer to this pair as the symbol $\sigma $ marked with $\delta $.
(b) We extend each pattern $a\in L\subseteq \Sigma ^{F}$ one or more
patterns $a'\in (\Sigma ')^{F}$ by superimposing new symbols over
each symbol of $a$. Call the new syntax $L'$. The SFT $X'$ defined
by $L'$ has the property that every pattern appearing in $X'$ consists
of a $\Delta $-pattern superimposed over a $\Sigma $-pattern, and
the $\Sigma $-pattern is admissible for $X$. 

Note that the map $\pi :X'\rightarrow \Sigma ^{Z^{2}}$ which erases
the superimposed layer of data maps $X'$ into a subsystem of $X$.
We say that $x\in X$ is represented in $X'$ if one can turn $x$
into a point of $X'$ by superimposing a suitable $\Delta $-pattern
over $x$; i.e., if $x=\pi (x')$ for some $x'\in X'$.

\section{\label{sec:base}Step 1: Constructing the Base}

In this section we construct a two-dimensional SFT $X$ whose symbols
are marked with the symbols $0,1$. The symbol $1$ may appear with
any density in points of $x$, but for each fixed $x\in X$ the density
of $1$'s will be extremely uniform.

\subsection{An almost periodic SFT}

Consider the substitution on the alphabet $\{\circ ,\bullet \}$ defined
the the rule\[
\begin{array}{cc}
  & \\
 \bullet  & \mapsto \end{array}
\, \left(\begin{array}{cc}
 \circ  & \bullet \\
 \bullet  & \circ \end{array}
\right)\qquad \qquad \begin{array}{cc}
  & \\
 \circ  & \mapsto \end{array}
\, \left(\begin{array}{cc}
 \circ  & \circ \\
 \bullet  & \circ \end{array}
\right)\]

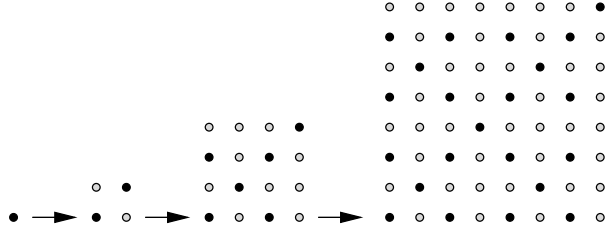
\begin{figure}

\input{2-net.pstex_t}

\caption{\label{fig:2-net} Three  iterations of the  substitution. }

\end{figure}

Let $W$ denote the dynamical system defined by these rules. See figure
\ref{fig:2-net}.

We say that a set $E\subseteq \mathbb{Z}^{2}$ is a \emph{$2$-net}
if $E=\cup _{n=1}^{\infty }I_{n}\times J_{n}$ where each $I_{n}$
and $J_{n}$ are translates of $2^{n}\mathbb{Z}$, the $I_{n}$'s
are pairwise disjoint, and the $J_{n}$'s are pairwise disjoint. We
refer to $I_{n}\times J_{n}$ as the \emph{$n$-th level} associated
with $E$. Note that if $u$ belongs to some level of $E$ then the
row and column to which $u$ belongs do not intersect any other level. 

\begin{prop}
Let $w\in W$ and\[
E(w)=\{u\in \mathbb{Z}^{2}\, :\, w(u)=\bullet \}\]
Then $E=E(w)$ is a $2$-net.
\end{prop}
\begin{proof}
For $n=0,1,2,\ldots $ let $a_{n}$ be the sequence of $F_{2^{n}}$
square patterns obtained by applying the substitution rule to the
initial symbol $a_{0}=\bullet $. It is sufficient to show that there
is a $2$-net $E=\cup _{n=1}^{\infty }I_{n}\times J_{n}$ such that
$\{u\in F_{2^{n}}\, :\, a_{n}(u)=\bullet \}=E\cap F_{2^{n}}$. To
verify this, one proves by induction that the above holds for \[
I_{n}=J_{n}=2^{n}\mathbb{Z}+2^{n-1}\qedhere \]
 
\end{proof}
We remark that the system $W$ supports a unique invariant probability
measure and as a measure preserving system this is an odometer, i.e.
is isomorphic to a zero-dimensional abelian group along with a free
minimal $\mathbb{Z}^{2}$ action generated by translation by two elements
of the group. 

This substitution rule has unique derivation, since one may check
that there is unique way to derive the central $6\times 6$ square
of the $8\times 8$ pattern in figure \ref{fig:2-net} from a $4\times 4$
pattern. 

Let $\widetilde{W}$ be the SFT associated to $W$ by theorem \ref{thm:subst}.
Then to each point in $\widetilde{W}$ there is associated, via a
one-block map, a $\{\circ ,\bullet \}$ pattern defining a 2-net.

\subsection{Marking the columns of $\widetilde{W}$}

We now superimpose another layer on top of $\widetilde{W}$. Begin
by superimposing the symbols $0,1$ on top of the $\widetilde{W}$
with the constraint that the symbols $0,1$ cannot be placed vertically
adjacent to each other. This forces each column in the resulting system
to be marked either entirely with $0$'s or entirely with $1$'s.

For a point $w\in \widetilde{W}$, the new coloring induces a $\{0,1\}$-coloring
of each level $I\times J$ in the decomposition given by the proposition.
This coloring is constant on the intersection of $I\times J$ with
columns; we now force it to be constant on the intersection of the
grid with rows. For this, superimpose two new symbols {}``$\longleftrightarrow $'',''$\Longleftrightarrow $''
on top of the existing ones. We think of $\longleftrightarrow $ as
transmitting a {}``$0$'' signal, and of $\Longleftrightarrow $
as transmitting a {}``$1$'' signal. The rules are that over a symbol
marked $\bullet $, the symbol $\longleftrightarrow $ appears always
together with the symbol $0$, and $\Longleftrightarrow $ appears
always together with the symbol $1$. We also require that $\longleftrightarrow $
and $\Longleftrightarrow $ cannot appear as horizontal neighbors,
so the arrow type is constant on rows.

Call the resulting system $X$ (it is of course an SFT) and let $x\in X$
be superimposed over a point $w\in \widetilde{W}_{0}$. Let $I\times J$
be some level of the $2$-net induced by $w$, and suppose that $w(u)$
is marked $0$ for some $u\in I\times J$. .Since it is also marked
$\bullet $, it bears the symbol $\longleftrightarrow $ (and not
$\Longleftrightarrow $); this forces the entire row to which $u$
belongs to be marked with $\longleftrightarrow $. Every other $v\in I\times J$
belonging to the same row is thus marked $\bullet $ and $\longleftrightarrow $,
and so it must be marked $0$. A similar analysis holds if $w(u)$
is marked $1$.

In short, the $0,1$-coloring of each grid $I\times J$ is constant
on rows and columns, and thus is completely constant. If $I_{n}\times J_{n}$
are the levels of the $2$-net induced by a point $x\in X$ then each
$I_{n}\times J_{n}$ determines a collection of columns which is $2^{n}$-periodic
in the horizontal direction, and all these columns bear the same symbol
$0$ or $1$. 

For $x\in X$, let $\delta (x)$ be the upper density of $1$'s in
$x$, i.e. \[
\delta (x)=\limsup _{n\rightarrow \infty }\frac{|\{u\in F_{n}:x(u)=1|\}}{|F_{n}|}\]
 where as usual $F_{n}=\{1,\ldots ,n\}^{2}$. If $I_{n}\times J_{n}$
are the levels of the $2$-net induced by $x$, then a simple calculation
shows that\[
\delta (x)=\sum _{n=1}^{\infty }\rho _{n}\cdot 2^{-n}\]
where $\rho _{n}$ is $0$ or $1$ according to the coloring $x$
induces on $I_{n}\times J_{n}$. Since the $I_{n}$'s and $J_{n}$'s
are pairwise disjoint the arrows transmitting information between
the points of each grid occupy different rows, and hence don't interact.
Therefore, we are free to color each level $0$ or $1$ independently
of the coloring of the other levels. Consequently, any sequence $\rho _{n}\in \{0,1\}$
may arise, so there are points $x\in X$ with $\delta (x)$ taking
on any value in the range $[0,1]$.

We will call a point in $X$ exceptional if it is superimposed over
an exceptional point of $\widetilde{W}$. For an exceptional point
$x\in X$ there are complementary half-spaces and/or quarter-spaces
such that the restriction of $x$ to each of them looks like a non-exceptional
point. Thus the above analysis applies to each of these regions separately.
This is not to say that we can glue admissible half- and quarter-spaces
together arbitrarily, and indeed for exceptional points the arrows
from different parts can interact; but this will not matter to us. 

Finally, we claim that $X$ has zero entropy. Indeed, $W$ has zero
entropy, and it is simple to check that if $a$ is a square pattern
admissible for $W$ then every extension of $a$ to a pattern $b$
admissible for $X$ is determined by the symbols of $b$ on the boundary
of the square. It follows that $X$ has entropy $0$.

\section{\label{sec:Pruning}Step 2: Pruning}

Let $h$ be a fixed right recursively enumerable number. Let $X$
be the system constructed in the previous section. Our goal in this
section is to construct an SFT $Y$ superimposed over $X$ which {}``kills''
points with density of $1$'s greater than $h$. More precisely, we
will want \[
\sup \{\delta (y)\, :\, y\in Y\}=h\]
 (here $\delta $ is the natural extension of $\delta $ from $X$
to $Y$) and that the supremum will be achieved.

\subsection{\label{sub:Boards}Boards}

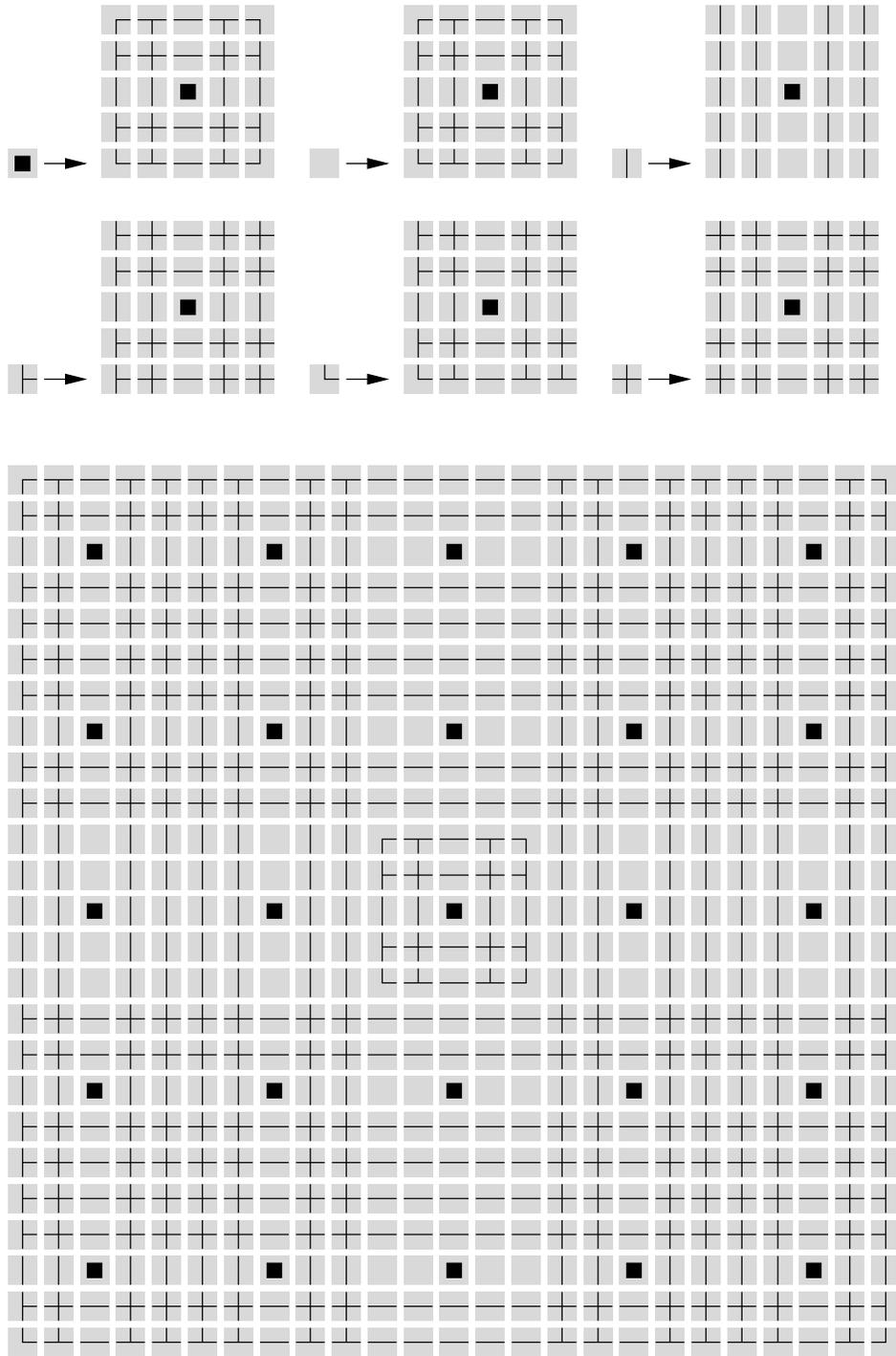
\begin{figure}

\input{5-board.pstex_t}

\caption{\label{fig:5-board}The substitution rules, up to rotation. The symbol $\Box$ is represented as an empty  square. The large 25x25 pattern is obtained by applying the substitution rules twice to $\blacksquare$ }

\end{figure}

We define a substitution system over the alphabet\[
\Sigma =\{|,-,\ulcorner ,\urcorner ,\llcorner ,\lrcorner ,\top ,\bot ,\vdash ,\dashv ,+,\blacksquare ,\Box \}\]
The substitution rules are described in figure \ref{fig:5-board}
together the symmetric rules obtained by rotating by multiples of
$90^{\circ }$. Let us denote by $b_{n}$ the $5^{n}\times 5^{n}$-pattern
obtained by applying the substitution rule $n$ times to the symbol
$\blacksquare $; see figure \ref{fig:5-board}. It is not hard to
show that $\blacksquare $ appears with period $5$ in every $b_{n}$.
Given $k$ and $n>k$, since $\blacksquare $ appear in $b_{n-k}$
with period $5$ we see that $b_{i}$ appears in $b_{n-k+i}$ with
period $5^{i}$, so $b_{k}$ appears in $b_{n}$ with period $5^{k}$. 

As can be seen from figure \ref{fig:5-board}, this substitution rule
produces patterns which induce certain grid-like shapes on $\mathbb{Z}^{2}$
. More precisely, define finite sets $I_{n}\subseteq \mathbb{N}$
inductively by $I_{1}=\{1,2,4,5\}$ and\[
I_{n+1}=I_{n}\cup (I_{n}+5^{n})\cup (I_{n}+3\cdot 5^{n})\cup (I_{n}+4\cdot 5^{n})\]
One sees by induction that $\min I_{n}=1$ and $\max I_{n}=5^{n}$,
so the union above is disjoint, and hence $|I_{n}|=4^{n}$ for each
$n$. Let \[
B_{n}=(I_{n}\times \{1,2,\ldots ,5^{n}\})\cup (\{1,2,\ldots ,5^{n}\}\times I_{n})\]
This is the set obtained by {}``filling in'' the rows and columns
between points of $I_{n}\times I_{n}$. The set $B_{n}$ is called
an $n$\emph{-board}. An infinite board is any set $B\subseteq \mathbb{Z}^{2}$
which is the limit of a sequence of translates $B_{n}+u_{n}$ for
some $u_{n}\in \mathbb{Z}^{d}$, where by limit we mean that $u\in B$
if and only if eventually $u\in B_{n}+u_{n}$. It is simple to check
that every infinite board $B$ has density zero, i.e. for every $\varepsilon >0$
there is an $N$ so $\frac{|B\cap (F_{N}+u)|}{|F_{N}|}<\varepsilon $
for every $u\in \mathbb{Z}^{2}$. This follows from the recursion
formula for $I_{n}$.

Consider the patterns $b_{n}|_{B_{n}}$. One shows by induction that
these patterns do not contain the symbols $\blacksquare ,\Box $.
Also, for $u\in B_{n}$ the points $v\in B_{n}$ which are adjacent
to $u$ -- i.e., which differ from $u$ by $\pm e_{1}$ or $\pm e_{2}$
-- are determined by $b_{n}(u)$ by interpreting the symbol $b_{n}(u)$
as a collection of lines pointing to the neighbors of $u$ in $B_{n}$.
Thus, $\bot $ indicates that there are neighbors left, right and
above the current symbol; $-$ indicates neighbors to the left and
right of it; etc. One can show that if $u\in F_{5^{n}}=\{1,\ldots ,5^{n}\}\times \{1,\ldots ,5^{n}\}$
and $b_{n}(u)\notin \{\blacksquare ,\Box \}$, then there is a unique
$k$ and translate $A$ of $B_{k}$ so that $u\in A$ and $b_{n}|_{A}$
is congruent to $b_{k}|_{B_{k}}$. In the large square in figure \ref{fig:5-board}
there are two boards visible; a $1$-board in the center, and a $2$-board
surrounding it. If we iterate the substitution one more step, each
$\blacksquare $ will turn into a $1$-board plus $\blacksquare $'s,
the $1$-board will turn into a $2$-board plus $\blacksquare $'s,
and the $2$-board will turn into a $3$-board, plus $\blacksquare $'s.

Let $R$ denote the dynamical system defined by these rules. From
the remarks above it follows that each $r\in R$ determines a pairwise
disjoint collection of boards, with $n$-boards appearing periodically
with period $5^{n}$; and if $u\in r$ and $r(u)\notin \{\blacksquare ,\Box \}$,
then $u$ belongs to one of these boards and the neighbors of $u$
in this board can be determined from symbols $r(u)$. By compactness,
there will exists points $r\in R$ and infinite boards $B$ so that
$r|_{B}$ is marked similarly to a finite board. Since infinite boards
cannot overlap and each occupies at least some quarter-space, there
can be at most four infinite boards in $r$, and since each has density
$0$, the density of points belonging to infinite boards in $r$ is
zero.

$R$ has unique derivation; indeed, the location of the corner tiles
determine the derivation of a point. We denote by $\widetilde{R}$
the SFT associated to $R$ by Mozes' theorem. We identify points in
$\widetilde{R}$ with the point in $R$ they are mapped to by the
given one-block map; in general this identification is many-to-one.

\subsection{Turing machines and their representations in SFTs}

A \emph{Turing machine} is an automaton with a finite number of internal
states which reads and writes data on a one-sided infinite array of
\emph{cells} indexed by $\mathbb{N}$, called the \emph{tape.} Each
cell contains one symbol from the data alphabet (so in our model the
input is an infinite sequence). The computation begins with the machine
located at the $0$-th (leftmost) cell and in a special initial state,
and the tape contains some data which is the input to the computation.
The state of the data tape along with the location and internal state
of the machine are called a \emph{configuration}; a configuration
uniquely determines all future configurations. The computation proceeds
in discrete time steps. At each iteration the machine is located at
some cell, reads the symbol written there and based on this data and
on its internal state, performs three actions: (a) it replaces the
current data symbol with a new one, (b) it moves one cell to the left
or to the right, and (c) it updates its internal state. The computation
may halt after a finite number of steps if the machine either moves
off the tape (steps left at cell $0$) or enters a designated state,
called the halting state. Barring these occurrences, the computation
continues forever. 

Although a very simple model, any algorithm written in a modern computer
programming language can be implemented as a Turing machine, and it
is generally accepted that any effective computation can be performed
by a Turing machine; this is Church's thesis. For background and basic
facts on this subject, see \cite{HU79}.

Let $X$ be the SFT constructed in section \ref{sec:base}, let $\widetilde{R}$
be the SFT described above and let $T$ be a Turing machine whose
data alphabet includes symbols $0,1$. We construct an SFT $Y_{T}$
superimposed over $X\times \widetilde{R}$ such that when a point
$y\in Y$ is superimposed over $(x,r)\in X\times \widetilde{R}$,
each board induced by $r$ has superimposed over it a pattern representing
the run of $T$ on the input given by the sequence of $0,1$'s appearing
in $x$ along the columns of the board. This construction, which we
describe next, is similar to the one used by Robinson in \cite{robinson71},
except that Robinson's machines always ran on an {}``empty'' input.

Let $\xi ,\rho $ be symbols in the alphabets of $X,\widetilde{R}$
respectively. We superimpose new symbols over $(\xi ,\rho )$ only
if $\rho $ represents a point in a board (i.e. $\rho \neq \Box ,\blacksquare $),
and the adjacency rules for the new symbols will only restrict pairs
of neighbors which belong to the same board (note that this can be
determined locally). Thus $(x,r)\in X\times \widetilde{R}$ will be
represented in $Y_{T}$ if and only if for each (finite or infinite)
board $B$ induced by $r$ there exists a locally admissible pattern
superimposed over $(x,r)|_{B}$. 

For a board $B_{n}+u$ let us call the points $I_{n}\times I_{n}+u$
the \emph{nodes} of the board. Note that $(\xi ,\rho )$ represents
a node if and only if $\rho \in \{\ulcorner ,\urcorner ,\llcorner ,\lrcorner ,\top ,\bot ,\vdash ,\dashv ,+\}$.
The data superimposed over a node will include a combination of data
symbol (from the machine's data alphabet) and possibly also a machine
state; this information may be represented by the alphabet $\Delta _{1}\cup (\Delta _{1}\times \Delta _{2})$
where the union is disjoint, $\Delta _{1}$ is the machine's data
alphabet and $\Delta _{2}$ its state space. 

Each row of nodes in a board is to represent a finite portion of the
configuration of the machine. More precisely, each node will contain
either a data symbol or a data symbol and a machine state; this is
called the cell's configuration. Suppose $x\in X$ and $r\in \widetilde{R}$
induces a board $B$. We can arrange things so that 

\begin{enumerate}
\item The data symbols in the nodes of the bottom row are the symbol $0$
or $1$ induced by $x$ on that node.
\item The node at the lower left corner of $B$ contains the initial state
of the machine, and no other node in the bottom row contains a machine
state.
\item Each row of nodes except the bottom one represents the configuration
obtained by iterating the computation one step from the configuration
given in the row below it. In particular, no row can appear admissibly
above a row containing a halting state.
\end{enumerate}
Properties (1) and (2) are easily implemented by restricting the types
of symbols which may be superimposed over $(\xi ,\rho )$ when $\rho \in \{\llcorner ,-,\perp ,\lrcorner \}$. 

Implementing (3) with local rules requires a little more effort. First,
note that in the course of the operation of a Turing machine $T$,
the configuration of a cell $i$ at a time $t>1$ is a function of
the configurations of the cells $i-1,i,i+1$ at time $t-1$; indeed
the data on the cell is determined by the configuration at $i$, and
the presence and state of the machine depend on the configurations
of the cells at $i-1,i+1$ (in case $i=0$, the dependence is on the
cells at $i,i+1$ only). We write $T(u,v,w)$ for the state of $i$
at time $t+1$ given that at time $t$ cells $i-1,i,i+1$ were in
states $u,v,w$ respectively (we allow $u=\textrm{"null"}$ in case
$i=0$). If we forget the geometry of the boards and imagine configurations
of the machine represented as sequences of cell configurations stacked
one on top of the other, this transition is {}``local'' and can
be enforced by a local rule that every pattern of the form $\begin{array}{ccc}
  & v' & \\
 u & v & w\end{array}
$ must satisfy $v'=T(u,v,w)$. 

However, when we represent cell configurations in nodes of a board
the transition from row to row is no longer local, since in a board
the nodes representing successive cells are spread out in space and
may be arbitrarily far apart. We can overcome this by using the rows
and columns between nodes (which belong to the board, and therefore
do not overlap for distinct boards) to {}``transmit information''.
In this way we can guarantee that the symbol superimposed over the
immediate neighbors of each node indicate the cell configuration at
each of the neighboring nodes. This can be implemented in a manner
similar to the way in which we synchronized the coloring of $2$-nets
in $X$ in section \ref{sec:base}. Briefly, over each grid point
marked $-$ we superimpose a pair of symbols $(u,v)$ where $u,v$
are node configurations. We require that each pair of horizontally
adjacent $-$'s are marked with the same pair, so all members of an
uninterrupted horizontal sequence of $-$'s carry the same pair. When
a pair $+-$ appear and $+$ has configuration $u$ we require that
over $-$ there is a pair $(u,v)$ for some $v$; and similarly for
pairs $\perp -$ and $\top -$. The symmetric condition is imposed
for $-+$, $-\perp $ and $-\top $. The result is that every uninterrupted
horizontal sequence of $-$'s carries the pair $(u,v)$ where $u$
is the configuration of the node at which the sequence ends on the
left, and $v$ the configuration of the node ending the sequence on
the right.

Next, over each symbol $|$ we superimpose a pair $(u,v,w)$, where
$u,v,w$ are cell configurations and $u$ or $w$ may also be {}``blank''.
As for $-$'s, we require that the marking is constant for each uninterrupted
vertical sequence of $|$'s. The markings are determined as follows.
If a $|$ is located immediately above a node with configuration $v$,
and the nodes to the left and right of that node have configurations
$u,w$ respectively, then $|$ carries $(u,v,w)$; $u$ or $w$ are
be {}``blank'' in the case there is no node to the left or right
of the node below $|$ (i.e. if it is at the edge of the board). Note
that by the previous discussion, $u,v,w$ may be determined by looking
at the immediate neighbors of the $|$. Thus the column of $|$'s
above each node represented the configuration of that node and its
neighbors. 

Finally, we require that when a node in state $v'$ appears vertically
above a $|$ marked $(u,v,w)$, then $v'=T(u,v,w)$. These conditions
can be seen to force property (3). 

We summarize this construction and its properties in the following
proposition:

\begin{prop}
\label{prop:YT-extension}\label{prop:YT-entropy}Given the systems
$X,\widetilde{R}$ from sections \ref{sec:base} and \ref{sub:Boards}
respectively, and given a Turing machine $T$, there exists an SFT
$Y_{T}$ superimposed over $X\times \widetilde{R}$ such that the
following are equivalent: 
\begin{enumerate}
\item $(x,r)\in X\times \widetilde{R}$ is represented in $Y_{T}$.
\item For each finite or infinite board $B$ induced by $r$ and containing
the symbol $\llcorner $, when $T$ is run on the sequence of $0,1$-s
induced by $x$ on the columns of $B$ the number of steps it runs
without halting is at least equal to the number of rows in $B$.
\end{enumerate}
Furthermore, $h(Y_{T})=0$.
\end{prop}
\begin{proof}
The equivalence follows from the discussion preceding the theorem.
We only note that if a board $B$ induced by $r$ does not contain
the symbol $\llcorner $ then it can always be extended, e.g. by a
pattern in which all rows are the same and contain only data. Note
that in general, there may be infinitely many ways to superimpose
a pattern over an infinite board which does not contain $\llcorner $.
Thus the projection from $Y_{T}$ into $X\times \widetilde{R}$ is
not an injection. 

It remains to check that $h(Y_{T})=0$. Given an $N\times N$ pattern
$a$ appearing in $X\times \widetilde{R}$, if $B_{n}+u$ is a board
induced by $\widetilde{R}$ and contained in $F_{N}$ then there is
a unique way to extend $a$ to a locally admissible $Y_{T}$ pattern.
This is true also for symbols in $a$ which do not lie in any board.
Given $\varepsilon >0$, a simple estimate shows that if $N$ is large
enough these points make up all but an $\varepsilon $-fraction of
the points in $F_{N}$, the remaining points coming from {}``boards''
which intersect the boundary of $F_{N}$ or infinite boards, all of
which have density tending to zero as $N\rightarrow \infty $. Hence
$a$ can be completed in at most $2^{\varepsilon (N)\cdot N^{2}}$
ways with $\varepsilon (N)\rightarrow 0$. It now follows that\[
N_{Y_{T}}(F_{n})\leq N_{X\times \widetilde{R}}(F_{n})\cdot 2^{\varepsilon (n)n^{2}}\]
therefore\[
h(Y_{T})\leq \frac{1}{n^{2}}\lim _{n\rightarrow \infty }N_{X\times \widetilde{R}}(F_{n})+\lim _{n\rightarrow \infty }\varepsilon (n)\leq h(X)+h(\widetilde{R})=0\]
as claimed.
\end{proof}

\subsection{\label{sub:Pruning}Pruning }

Our aim now is to find a Turing machine $T$ so that $(x,r)\in X\times \widetilde{R}$
is represented in $Y_{T}$ if and only if $\delta (x)\leq h$. 

Recall that this machine $T$ will receive as its input sequences
of $0,1$'s induced by points $x\in X$ on translates of $I_{n}$.
Write $I=\cup _{n=1}^{\infty }I_{n}$, and enumerate the elements
of $I$ as $I=\{i_{1}<i_{2}<\ldots \}$, where $i_{1}=1$. Note that
the first $4^{n}$ elements of this sequence are precisely the elements
of $I_{n}$; this follows easily from the recursion relation defining
the $I_{n}$'s. If $(x,r)\in X\times \widetilde{R}$ and $B_{n}+u$
is an $n$-board induced by $r$, then the $0,1$-coloring induced
by $x$ on $I_{n}\times \{1\}+u$ is the sequence $(x_{j})_{j=1}^{4^{n}}$
such that $x_{j}$ is the symbol $0$ or $1$ appearing on the $i_{j}$-th
column in $T^{u}x$. It follows that for any $k\leq n$, the first
$4^{k}$ symbols of this sequence correspond to a pattern induced
by $x$ on some translate of $I_{k}$. 

\begin{lem}
\label{lem:equi-dist-B_n}There is a sequence of finite sets $M_{n}\subseteq I_{m\cdot 2^{n^{2}}}$
such that $\{i_{m}\, :\, m\in M_{n}\}$ is a complete set of residue
classes modulo $2^{n}$, i.e. for every $0\leq j<2^{n}$ there exists
a unique $m\in M_{n}$ such that $i_{m}\equiv j\bmod 2^{n}$.
\end{lem}
\begin{proof}
By the recursion formula for $I_{k}$ given section \ref{sub:Boards}
and the fact that $I_{k}$ is an increasing sequence, for any $k\leq q$
we have\[
I_{k}+5^{q}\subseteq I_{q}+5^{q}\subseteq I_{q+1}\]
 In particular, since $1\in I_{1}$, we may show by induction that
for any $r$ and $t$, \[
1+5^{q}+5^{2q}+\ldots +5^{tq}\subseteq I_{q+2q\ldots +tq+1}\]
Given $n$, since $\gcd (2^{n},5)=1$ we may choose $q\leq 2^{n}$
so that $5^{q}\equiv 1\bmod 2^{n}$. Since the set \[
\{1+5^{q}+\ldots +5^{tq}\}_{t=1}^{2^{n}}\]
 is a complete set of residues modulo $2^{n}$ and is contained in
$I_{2^{3n}}$; the existence of $M_{n}$ follows.
\end{proof}
It is clearly possible to compute a sequence of sets $M_{n}\subseteq I_{2^{3n}}$
with the above properties. The proof above gives an algorithm for
going so, since the identity $5^{q}=1\bmod 2^{n}$ is solved by $q=\phi (2^{n})$
(here $\phi $ is Euler's function).

Let $r(n)$ be a computable sequence and $h\leq r(n)\rightarrow h$.
We can now describe our algorithm:

\begin{algorithm}
\label{alg:main}Input: $(x_{n})_{n\in \mathbb{Z}}\in \{0,1\}^{\mathbb{Z}}$.

For $N=1,2,3\ldots $ do
\begin{enumerate}
\item Calculate $r(N)$.
\item Calculate the relative frequency $\delta _{N}$ of $1$'s in the sequence
$(x_{m}\, :\, m\in M_{N})$, i.e.\[
\delta _{N}=\frac{1}{2^{N}}\#\{m\in M_{N}\, :\, x_{m}=1\}\]

\item If $\delta _{N}>r(N)+2^{-N}$ then halt.
\end{enumerate}
\end{algorithm}
\begin{prop}
\label{pro:algo-analysis}Let $x\in X$ and let $(x_{n})_{n=1}^{\infty }$
be the $0,1$-valued sequence with $x_{n}$ equal to the color of
the $i_{n}$-th column of $x$. Then the algorithm \ref{alg:main}
halts on the input $(x_{n})$ if and only if $\delta (x)>h$, and
if it halts the number of steps it runs before halting depends only
on $\delta (x)$ (not on $x$).
\end{prop}
\begin{proof}
It suffices to show that $\delta (x)-2^{-N}\leq \delta _{N}\leq \delta (x)+2^{-N}$
for every $N$. Indeed, if $\delta (x)=h+\varepsilon $ for some $\varepsilon >0$
then $\delta _{N}\geq \delta (x)-2^{-N}$ implies that $\delta _{N}>h-\varepsilon /2$
for large enough $N$, and since $r(N)\rightarrow h$ for large enough
$N$ we will have $\delta _{N}>r(N)+2^{-N}$ and the algorithm will
halt. On the other hand if $\delta (x)\leq h$ then $\delta _{N}\leq \delta (x)+2^{-N}$
implies that $\delta _{N}\leq h+2^{-N}\leq r(N)+2^{-N}$, so the algorithm
will run forever.

Fix $N\geq 1$ and let $E=\cup _{n=1}^{\infty }U_{n}\times V_{n}$
be the $2$-net induced by $x$. Note that \[
\delta (x)=\sum _{n=1}^{\infty }\rho _{n}2^{-n}\]
where $\rho _{n}\in \{0,1\}$ is the symbol induced by $x$ on the
grid $U_{n}\times V_{n}$.

Note that $|M_{n}|=2^{n}$. Let $J_{n}=\{j\in M_{N}\, :\, i_{j}\in U_{n}\}$.
Since $\{i_{j}\, :\, j\in M_{N}\}$ is a complete set of residues
modulo $2^{N}$, for each $n\leq N$ we have \[
|J_{n}|=2^{N-n}\]
 and since the $U_{n}$'s are pairwise disjoint so are the $J_{n}$'s,
so \[
|M_{N}\setminus \bigcup _{n=1}^{N}J_{n}|=1\]
Let $M_{N}\setminus \cup _{n=1}^{N}J_{n}=\{i\}$ and let $\rho '\in \{0,1\}$
be the symbol induced on the $i$-th column of $x$. Then\begin{eqnarray*}
\delta _{N} & = & \frac{1}{2^{N}}\sum _{j\in M_{N}}1_{\{x_{j}=1\}}\\
 & = & \frac{1}{2^{N}}(1_{\{x_{i}=1\}}+\sum _{n=1}^{N}\sum _{j\in J_{n}}1_{\{x_{j}=1\}})\\
 & = & \frac{\rho '}{2^{N}}+\sum _{n=1}^{N}\frac{\rho _{n}|J(n)|}{2^{N}}\\
 & = & \frac{\rho '}{2^{N}}+\sum _{n=1}^{N}\rho _{n}2^{-n}
\end{eqnarray*}
The desired inequality follows.

Regarding the number of steps the algorithm runs before halting, this
depends only on $N$ and $\delta (x)$.
\end{proof}
Let $T$ be a Turing machine implementing this algorithm and whose
input is the sequence of $0$'s and $1'$s which is the input to the
algorithm. We make one important assumption about the implementation,
namely that there are integers $t_{N}$ such that the machine performs
the first $N$ iterations of the loop in at most $t_{N}$ steps (or
halts before that), independent of the input. Such an implementation
does not present any difficulty. Another thing to note is that as
we have defined it, the entire tape is taken up by input data. In
order to provide the machine with space to store its intermediate
calculations one can allow it to superimpose another layer of symbols
over the input alphabet. Formally, this can be done by setting the
machines alphabet to be $\{0,1\}\times \{0',1'\}$, with the input
represented by the first coordinate and the machines modifying the
second coordinate as it pleases.

Let $Y=Y_{T}$; this is the system $Y$ whose construction was the
goal of the second step in the outline given in section \ref{sec:Outline}.

\begin{prop}
If $(x,r)\in X\times \widetilde{R}$ then $(x,r)$ is represented
in $Y_{T}$ if and only if $\delta (x)\leq h$. 
\end{prop}
\begin{proof}
By \ref{prop:YT-extension} it suffices to show that the condition
$\delta (x)\leq h$ is equivalent to the fact that for any finite
or infinite board $B$ induced by $r$ representing an $N\times N$
grid ($N\in \mathbb{N}\cup \{\infty \}$), if $r|_{B}$ contains the
symbol $\llcorner $ then the algorithm does not halt after $N$ steps
when run on the input induced by $x$ on the columns of $B$. The
proposition now follows easily from proposition \ref{pro:algo-analysis}
and the fact that $r$ induces boards of arbitrarily large size. 
\end{proof}
Finally, we note that the topological entropy of $Y=Y_{T}$ is zero
by proposition \ref{prop:YT-entropy}.

\section{\label{sec:calculating-entropy}Step 3: Adding and calculating entropy}

Let $Y$ be the system constructed in the previous section. Let $Z$
be the SFT superimposed over $Y$ by adding one of the symbols $\alpha ,\beta $
over each occurrence of the symbol $1$. We place no other restrictions
on the configurations of $\alpha ,\beta $'s which may appear. In
this section we estimate the entropy of $Z$ and show that it is indeed
equal to $h=\sup \{\delta (y)\, :\, y\in Y\}$.

Write $F_{n}=\{1,\ldots ,n\}^{2}$ again and for $y\in Y$ denote
\[
f_{n}(y)=\frac{1}{|F_{n}|}\#\{u\in \mathbb{Z}^{2}\, :\, y(u)\textrm{ is marked }"1"\}\]
 so $\delta (y)=\limsup _{n\rightarrow \infty }f_{n}(y)$. Since $\delta (y)\leq h$
for every $y\in Y$ there is a sequence $\varepsilon _{n}\rightarrow 0$
such that\[
\sup _{y\in Y}f_{n}(y)<h+\varepsilon _{n}\]
(such a sequence $\varepsilon _{n}$ exists by general considerations,
but in our case by the proof of proposition \ref{pro:algo-analysis}
one can choose $\varepsilon _{n}=2^{-n+1}$).

We now estimate the number of patterns induced by $Z$ on the box
$F_{n}=\{1,\ldots ,n\}^{2}$. For each pattern induced on $F_{n}$
by $y\in Y$, the number of ways to superimpose the symbols $\alpha ,b$
and get an admissible pattern for $Z$ is $2^{f_{n}(y)|F_{n}|}=2^{f_{n}(y)n^{2}}$.
Summing over all patterns induced on $F_{n}$ by $Y$ and using the
fact that $f_{n}(y)\leq h+\varepsilon _{n}$ we have \[
N_{Z}(F_{n})\leq N_{Y}(F_{n})\cdot 2^{f_{n}(y)n^{2}}\leq N_{Y}(F_{n})\cdot 2^{n^{2}(h+\varepsilon _{n})}\]
 so\[
\limsup _{n\rightarrow \infty }\frac{1}{|F_{n}|}\log N_{Z}(F_{n})\leq \limsup _{n\rightarrow \infty }\frac{1}{|F_{n}|}\log N_{Y}(F_{n})+\limsup _{n\rightarrow \infty }(h+\varepsilon _{n})=h\]
 because $h(Y)=\limsup \frac{1}{|F_{n}|}\log N_{Y}(F_{n})=0$.

On the other hand, if $y_{n}\in Y$ satisfy $f_{n}(y_{n})\rightarrow h$
then clearly the number of ways to extend $y_{n}|_{F_{n}}$ to a pattern
in $Z$ is $2^{f_{n}(y_{n})|F_{n}|}$ and so \[
\liminf _{n\rightarrow \infty }\frac{1}{|F_{n}|}\log N_{Z}(F_{n})\geq \limsup _{n\rightarrow \infty }f_{n}(y_{n})=h\]
 The entropy estimate $h(Z)=h$ follows. 

This completes the proof of theorem \ref{thm:main}.

\section{\label{sec:open-problems}Concluding remarks }

Many questions remain about the relation between the dynamics SFTs
and their entropies. Let us take a closer look at the system $Z$
constructed above. We can write $Z$ as a disjoint union $Z=\cup _{0\leq r\leq h}Z_{r}$
where $Z_{r}$ is the (nonempty) set of points $z\in Z$ with $\delta (z)=r$;
each $Z_{r}$ is a closed shift-invariant set, so every orbit closure
in $Z$ lies in some $Z_{r}$. Hence $Z$ is not transitive. $Z$
also does not have periodic points, since it factors onto the infinite
uniquely ergodic system $W$.

We remark that if $h$ is computable instead of merely right recursively
enumerable, then one can modify algorithm \ref{alg:main} so as to
also kill points whose density of $1$'s is less than $h$ (computability
implies both right and left recursive enumerability). For this algorithm
the resulting system is essentially the system $Z_{h}$ above. However,
it is still not transitive, since there are many ways to extend an
infinite board which does not contain a bottom row; this does not
affect entropy, since infinite boards have density zero, but means
that $Z_{h}$ has a transient part. 

\begin{problem}
Is every right recursively enumerable number $h$ the entropy of a
transitive SFT?
\end{problem}
Conversely, we have seen that the entropy of irreducible SFTs is computable.
This raises the following:

\begin{problem}
What is the class of entropies of multidimensional irreducible SFTs?
\end{problem}
Another mechanism which may be related to entropy is the presence
of periodic points. For a two-dimensional SFT $X$ denote by $P_{n}$
the number of $n\times n$ patterns which can be repeated to produce
an admissible tiling of the lattice with period $k$ in both directions.
Clearly $P_{n}$ is computable, and in certain situations one can
show that $\lim \frac{1}{n^{2}}\log P_{n}\rightarrow h$ and $\frac{1}{n^{2}}\log P_{n}\leq h+\varepsilon _{n}$
for a sequence $\varepsilon _{n}$ which decays to $0$ at a computable
rate. This implies that the entropy is computable, because for $\widetilde{N}_{n}$
as in section \ref{sec:Entropies-of-SFTs-are-co-RE} we have $h\in (\frac{1}{n}\log P_{n}-\varepsilon _{n},\frac{1}{n}\log \widetilde{N}_{n})$,
and so given $n$ we can examine the difference $\frac{1}{n}\log \widetilde{N}_{n}-(\frac{1}{n}\log P_{n}-\varepsilon _{n})$
for $n=1,2,3\ldots $, stop the first time it is less than $\frac{1}{n}$,
and give $\frac{1}{n}\log P_{n}-\varepsilon _{n}$ as our estimate.

Friedland \cite{friedland97} used this observation to deduce that
if the syntax of an SFT enjoys a certain spacial symmetry then the
entropy is computable. We note also that irreducible SFTs in two dimensions
have dense periodic points, but whether this is so in higher dimensions
seems to be open \cite{W94}. 

\begin{problem}
Do dense periodic points for an SFT imply that the entropy is computable?
\end{problem}
Finally, we repeat here an old question which we mentioned in the
introduction:

\begin{problem}
Is every sofic shift a factor of an SFT with the same entropy? 
\end{problem}
\bibliographystyle{plain}
\bibliography{entropies-of-SFTs}

\end{document}

%% file: 2-net.pstex_t
\begin{picture}(0,0)%
\includegraphics{2-net.pstex}%
\end{picture}%
\setlength{\unitlength}{4144sp}%
\begingroup\makeatletter\ifx\SetFigFont\undefined%
\gdef\SetFigFont#1#2#3#4#5{%
  \reset@font\fontsize{#1}{#2pt}%
  \fontfamily{#3}\fontseries{#4}\fontshape{#5}%
  \selectfont}%
\fi\endgroup%
\begin{picture}(3572,1333)(4537,-4738)
\end{picture}%

%% file: 5-board.pstex_t
\begin{picture}(0,0)%
\includegraphics{5-board.pstex}%
\end{picture}%
\setlength{\unitlength}{4144sp}%
\begingroup\makeatletter\ifx\SetFigFont\undefined%
\gdef\SetFigFont#1#2#3#4#5{%
  \reset@font\fontsize{#1}{#2pt}%
  \fontfamily{#3}\fontseries{#4}\fontshape{#5}%
  \selectfont}%
\fi\endgroup%
\begin{picture}(5582,8473)(4230,-11402)
\end{picture}%